\documentclass[twoside]{oss-conf-eng}     
\usepackage{lmodern}
\usepackage{cmap}
\usepackage[english,slovak]{babel}             
\usepackage[T1]{fontenc}
\usepackage[utf8]{inputenc}
\usepackage{amsmath}
\usepackage{amscd,amssymb,amsfonts}
\usepackage[sort,nocompress]{cite} 
\usepackage{graphicx}                          
\usepackage{fancyhdr}
\usepackage[pdftex,unicode,bookmarks=false]{hyperref}  
\usepackage{listings}
\usepackage{eurosym}

\usepackage{tikz}
\usepackage{floatflt}

\begin{document}
\def\keywordsnameS{Kľúčové slová}     

\setcounter{page}{1}                            
\def\volumeDOI{OSSConf 2012:\  }                
\def\konfera{Konferencia OSSConf 2012}

\pagestyle{fancy}
\fancyfoot{}
\fancyhead[LE,RO]{\thepage}
\fancyhead[LO]{\nouppercase{\leftmark: \rightmark}}
\fancyhead[RE]{\nouppercase{\konfera}}

\title[Krivka konštantných uhlov viditeľnosti]{Krivka konštantných uhlov viditeľnosti\\ pre niektoré konvexné množiny v rovine}

\titleA[Equal aperture angles for some convex plane sets]{EQUAL APERTURE ANGLES CURVE\\ FOR SOME CONVEX SETS IN THE PLANE}

%

\author[M. Kaukič]{Michal Kaukič}{Mgr.}{CSc.}{}
\address{Department of Mathematical Methods, Faculty of Management Science and Informatics, 
University of Žilina, Univerzitná~8215/1, 010~26~Žilina, Slovak Republic}
\emailh{Michal.Kaukic@fri.uniza.sk}
\urladdressh{http://feelmath.info}


\keywordsA{Aperture angle, planar convex sets, equal aperture angle.}

\selectlanguage{slovak}

\begin{abstract}
Pre danú konvexnú množinu $K$ v rovine (môže byť aj neohraničená) môžeme zostrojiť krivku $C$, pre ktorú je uhol viditeľnosti tejto množiny konštantný, teda nadobúda predpísanú hodnotu. V článku uvádzame implicitný vzorec na výpočet $C$, zamýšľame sa nad praktickými spôsobmi výpočtu tejto krivky a uvádzame niekoľko jednoduchých príkladov, kedy sa dajú krivky rovnakých uhlov viditeľnosti určiť explicitne. Nakoniec je naznačené, aké sú ďalšie možné smery výskumu v tejto oblasti. V článku je v podstatnej miere použitý otvorený softvér (Sage, Pylab, IPython,\dots).
\end{abstract}

\selectlanguage{english}

\begin{abstractA}
For given convex set $K$ in the plane (not necessarily bounded), we can construct the curve $C$ for which the visibility (aperture) angle of this set has the same, prescribed value. We give the implicit formula for $C$, discuss some issues concerning practical computations of $C$ and bring several simple examples, when the equal visibility angle curves can be effectively computed explicitly. We conclude with some remarks about possible directions for further research in this area. Extensive use of Open Source software (Sage, Pylab, IPython,\dots) is a key feature of this article.
\end{abstractA}

\selectlanguage{english}
\newtheorem{theorem}{Theorem}[section]
\newtheorem{corollary}[theorem]{Corollary}
\newtheorem{lemma}[theorem]{Lemma}
\newtheorem{exmple}{Example}
\newtheorem{defn}{Definition}
\newtheorem{proposition}[theorem]{Proposition}
\newtheorem{conjecture}[theorem]{Conjecture}
\newtheorem{rmrk}[theorem]{Remark}
\newenvironment{definition}{\begin{defn}\normalfont}{\end{defn}}
\newenvironment{remark}{\begin{rmrk}\normalfont}{\end{rmrk}}
\newenvironment{example}{\begin{exmple}\normalfont}{\end{exmple}}
\newtheorem*{remarque}{Remark}
\newcommand\bbreak{\allowdisplaybreaks}
\newcommand\dd{\mathop{\rm d\!}\nolimits}
\newcommand\sgn{\mathop{\rm sgn}\nolimits}

\maketitle


\section{Introductory remarks}

\begin{definition}
Let $Q$ be the convex set in the plane. The aperture (or visibility) angle $\varphi(X, Q)$ of a point $X \not\in Q$ with respect to the set $Q$ is the angle of the smallest cone with apex $X$ that contains $Q$ (see e.g., \cite{angleapprox}, \cite{angleopt}). In this paper we will assume that the set $Q$ is closed, but not necessarily bounded.

For given angle $\alpha,\ 0 \le \alpha < \pi$ and the convex planar set $Q$ we will denote by $C(\alpha, Q)$ the set of all points $X$ for which the aperture angle $\varphi(X, Q)$ has the value $\alpha$.
\end{definition}
\begin{example}
Given the wedge $W$ with internal angle $\vartheta$, for fixed angle $\alpha$, three cases can occur: 
\begin{enumerate}
\item{for $\alpha < \vartheta$ the set $C(\alpha, W)$ is empty,}
\item{if $\alpha = \vartheta$ then $C(\alpha, W)$  is the wedge $V$ ($W$ rotated by angle $\pi$ about its apex), see Figure \ref{wedgea}, }
\item{for $\alpha > \vartheta$ the set $C(\alpha, W)$ is  the angle $V$ with value $2\alpha -\vartheta$, opposite to wedge $W$ and having common bisector with angle $\vartheta$, see Figure \ref{wb}.}
\end{enumerate}
\label{ex1}
\end{example}
\begin{figure}[ht]
\begin{minipage}[t]{5.5cm}
\begin{tikzpicture}[scale=0.8]
  \draw[fill=yellow!40] (15:2.5cm) -- (0,0) -- (60:2.5cm);
  \draw[fill=yellow!80] (0,0) -- (15:.75cm) arc (15:60:.75cm);
  \draw(40:0.47cm) node {$\vartheta$};  
  \draw(38:1.7cm) node {\large{$W$}};  

  \draw[fill=green!30] (15:-2.5cm) -- (0,0) -- (60:-2.5cm);
  \draw[fill=green!80] (0,0) -- (15:-.75cm) arc (15:60:-.75cm);
  \draw(35:-0.49cm) node {$\alpha$};  
  \draw(38:-1.7cm) node {\large{$V$}}; 
  \begin{scope}[thick]
    \draw (60:-2.5cm) -- (60:2.5cm);
    \draw (15:-2.5cm) -- (15:2.5cm);
  \end{scope}
\end{tikzpicture}
\caption{$C(\alpha, W) = V $ for $\alpha=\vartheta$}
\label{wedgea}
\end{minipage}
\hspace{1cm}
\begin{minipage}[t]{5.5cm}
\begin{tikzpicture}[scale=0.8]
  \draw[fill=yellow!40] (15:2.5cm) -- (0,0) -- (60:2.5cm);
  \draw[fill=yellow!80] (0,0) -- (15:.75cm) arc (15:60:.75cm) -- cycle;
  \draw[thick] (15:2.5cm) -- (0,0) -- (60:2.5cm);
  \draw(40:0.47cm) node {$\vartheta$};  
  \draw(38:1.7cm) node {\large{$W$}};  
  
  \draw(0,0) -- (-5:-1.65cm) arc (-5:80:-1.65cm);
  \draw[thick](-5:-2cm) -- (0,0) -- (80:-2cm);
  \draw(30:-0.93cm) node {\small{$2\alpha-\vartheta$}};  
  \draw(-5:-1.3cm) node[above] {\large{$V$}};
  \draw(80:-1.3cm) node[right] {\large{$V$}};
\end{tikzpicture}
\caption{$C(\alpha, W) = V $ for $\alpha>\vartheta$}
\label{wb}
\end{minipage}
\end{figure}

In the Example \ref{ex1} the convex set $W$ was unbounded and not strictly convex. The set $C(\alpha,W)$ was two-dimensional in one ``singular'' case. Next, we will consider unbounded, strictly convex set.

\begin{example}
Let us define the convex set $Q$ as $Q =\{ (x,y)\in {\mathbb{R}}^2: \ y \ge x^2 \}$, i.e., $Q$ is the parabola $y=x^2$ together with its interior points.

In this case, it is not difficult to compute the set $C(\alpha,Q)$ explicitly.
For arbitrary exterior point $A=(x_1,y_1) \not\in Q$ we can construct exactly two tangent lines to $Q$ passing through $A$. The slopes $k_1, k_2$ of that tangent lines are
\begin{equation}
k_1=\tan\varphi_1=2\left({x_1+\sqrt{x_1^2-y_1}}\right)\quad  
  k_2=\tan\varphi_2=2\left({x_1-\sqrt{x_1^2-y_1}}\right),
\label{kacka}  
\end{equation}  
where $\varphi_1,\varphi_2$ are the slope angles of two tangent lines. For given angle $\alpha$ we seek the set of all points $A$ such that
\begin{equation}
\tan(\varphi_1 - \varphi_2)=
\frac{\tan \varphi_1 - \tan \varphi_2}{1+\tan \varphi_1 \tan \varphi_2}=
 \frac{k_1-k_2}{1+k_1 k_2}=\tan\alpha=K.
\label{tangensy}
\end{equation}
If $k_1 k_2 \ne -1$, we can substitute for $k_1, k_2$ from (\ref{kacka}) and the solutions of resulting quadratic equation give us the explicit formula for the curve $C(\alpha, Q)$
\begin{equation}
y_p(x)=\frac{-K^2-2 \pm\, 2\sqrt{4 K^2 x^2 +K^2+1}}{4 K^2},
\label{parabcurve}
\end{equation}
where the plus sign is valid for angles $\alpha>\pi/2$ and the sign minus for sharp angles $\alpha$. For $\alpha=\pi/4$ this formula is simply
$y_p(x)=-(2\sqrt{4 x^2 +2}+3)/4$ and this is the curve with asymptotes $y=\pm\, x$. For arbitrary $K\ne 0$, the asymptotes of curve given by formula (\ref{parabcurve}) are $y=\pm\, x/K$.

For the singular case $k_1 k_2 = -1$, i.e., for $\alpha=\pi/2$ we obtain immediately $y_p(x)=-1/4$, which is the directrix of parabola $y=x^2$. From the formula (\ref{parabcurve}) we can see that for angle $\alpha<\pi/2$ the function $y_p(x)$ is concave and for  $\alpha>\pi/2$ we obtain convex function $y_p(x)$.

On the Figure \ref{parab} we have plotted the graphs of the functions $y_p(x)$ for several angles $\alpha$. From the bottom to top curve, the angles are: 
$\pi/4,\ \pi/3.5,\ \pi/3,\ \pi/2.5,$ $\pi/2,\ \pi/1.5,\ \pi/1.2$.  
\label{parabex}
\end{example}
\begin{floatingfigure}[p]{6.9cm}
       \includegraphics[width=6.5cm]{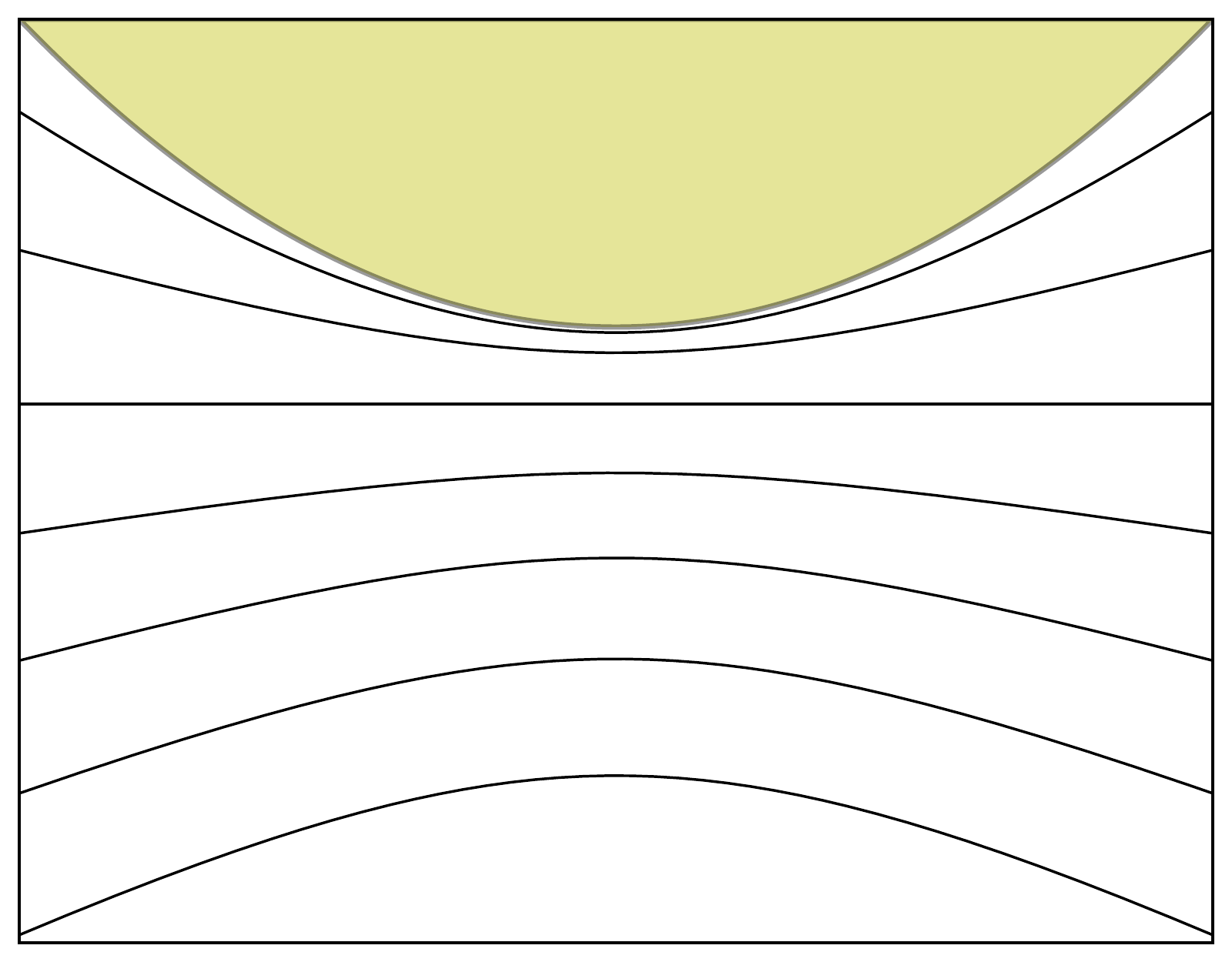}
       \caption{Functions $y_p(x)$ for parabolic region}
       \label{parab}
\end{floatingfigure}

Now, the set $C(\alpha,Q)$ is nonempty for all angles $0\le \alpha <\pi$. For compact (in $\mathbb{R}^2$ -- closed and bounded) $Q$ the set  $C(\alpha,Q)$ (if nonempty) forms the closed curve, which we will call {\bf constant visibility angle curve}. This curve have been investigated mainly for the special case of convex polygonal sets $Q$ (see \cite{teichmann}, \cite{angleapprox}, \cite{angleopt}). The papers cited are concerned with optimization problems (maximization) for visibility angle. In the paper \cite{anglethreed} the authors solved three-dimensional problem of maximiza\-tion of visibility angle for convex polyhedron, viewed from given line segment. Good sources of information about related problems are the Pirzadeh's Master thesis \cite{cgcalipers} and the article~\cite{gpcalipers}.    

\section{Constant visibility angle curves for convex functions}
We can, in principle, compute the  visibility angle curve for arbitrary convex function. Let us assume that the function $y=f(x)$ is strictly convex on finite segment $I=[a,b]$. Further, we assume that $f(x)$ has continuous first and second derivative on $I$. 

Given the (arbitrary, but fixed) point $\tilde{x} \in I$ and the angle $\alpha$, the problem is to find another point $x \in I$ such that (cf. equation (\ref{tangensy}))
\begin{equation}
\frac{f^\prime(\tilde{x})-f^\prime(x)}{1+f^\prime(\tilde{x}) f^\prime(x)}=\tan \alpha =K
\mbox{ or } f^\prime(x)=\frac{f^\prime(\tilde{x})-K}{1+K f^\prime(\tilde{x})}=M.
\end{equation} 
For given $\tilde{x}$, the constant $M$ is known, therefore we have to solve simple equation $F(x)=f^\prime(x)-M=0$. For some simple functions we can do that analytically but, in general, it is necessary to use one of suitable numerical methods. For example, using Newton method we get the iterational sequence:
$$\mbox{choose }\ x_0 \in I,\quad
x_{n+1}= x_n - \frac{f^\prime(x_n)-M}{f^{\prime\prime}(x_n)}\quad \mbox{ for } n=0,1,\dots$$ 
We will not analyze the convergence assumptions here. It will be the subject of further research. In the next section we give two nontrivial examples of computing the constant visibility curve for  smooth, bounded convex sets.

\section{Constant visibility angle curves for sine  and ellipse}

\begin{example}
Let $Q$ be the convex region, described by 
$$Q=\{ (x,y) \in {\rm R}^2:\ 0\le x \le \pi, \ -\sin x \le y \le 0\}.$$
We will take into account only visibility ``from below'', i.e., from points with negative $y$ coordinate.
\begin{floatingfigure}[p]{6.2cm}
       \includegraphics[width=5.7cm]{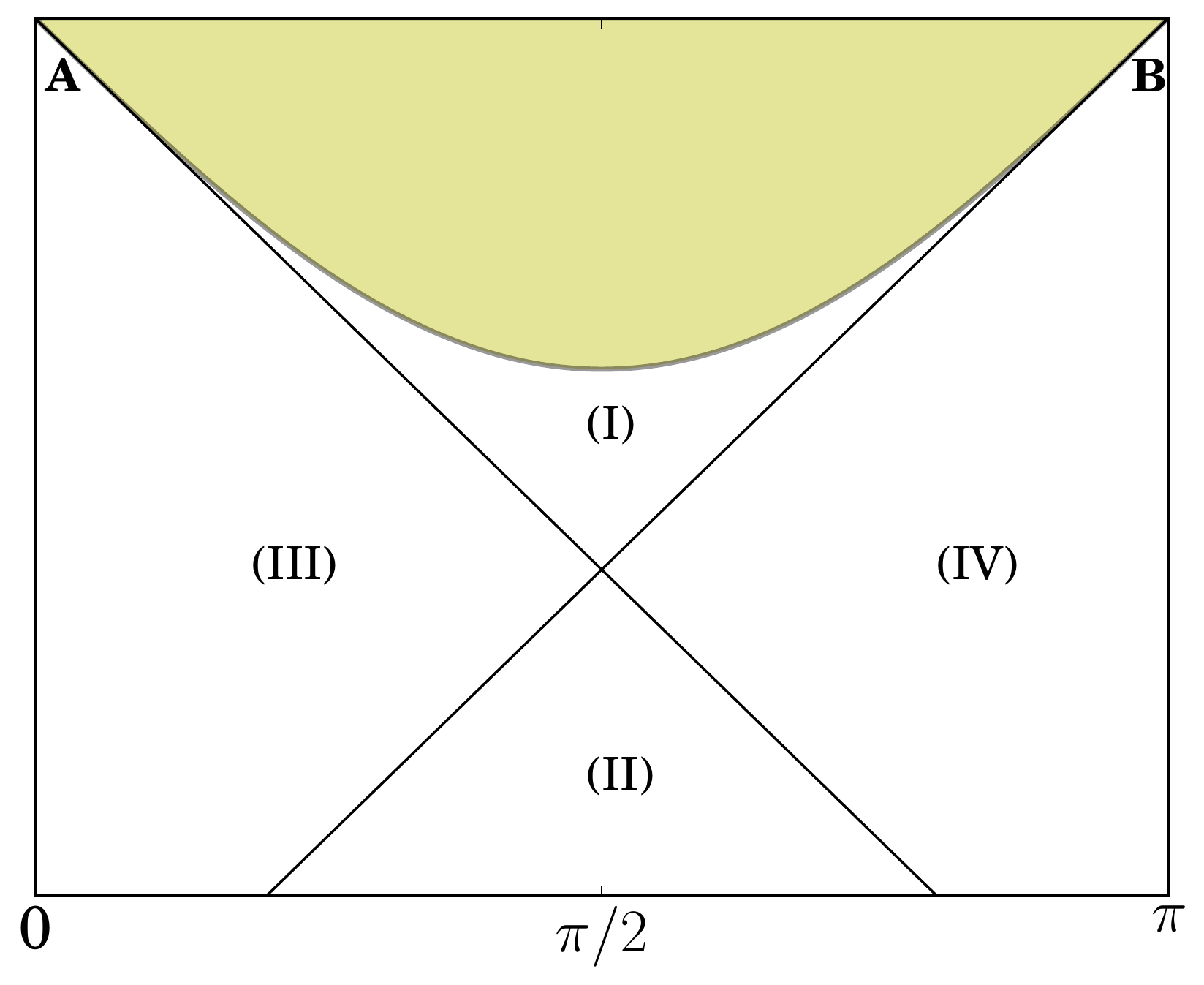}
       \caption{Regions for computing $C(\alpha,Q)$}
       \label{sindom}
\end{floatingfigure}
On the Figure \ref{sindom} we can see four main regions in the lower halplane. In each of them the computation of $C(\alpha,Q)$ is a bit different. From the symmetry of $Q$ it follows that it is sufficient to investigate only parts of regions with $x \le \pi/2$. In the region (I) there are only the points of curves $C(\alpha,Q)$ for $\alpha\ge \pi/2$. In the regions (III) resp. (IV) the visibility angles curve always contain points {\bf A} resp. {\bf B}. For the points in region (II) we have $\alpha \le \pi/2$, the set $Q$ is visible only as the segment $AB$, thus the curves $C(\alpha,Q)$ are circular arcs.
\end{example}

In this example, the curves $C(\alpha,Q)$ can be computed analytically, but to compute them ``by hand'' as in Example (\ref{parabex}) is very tedious and error-prone.
For symbolic computations, we utilized Open Source system Sage \cite{sage}. For numerical experiments the programming language Python~\cite{python}, with user friendly shell IPython~\cite{ipython} and modules Numpy~\cite{numpy}, Scipy~\cite{scipy} were used. Finally, nearly all graphics in this article was generated by   Python module Matplotlib~\cite{matplotlib}.   

Programs used for computations and generation of graphics can be downloaded from \url{http://feelmath.info/images/constangle/constangle.zip}. For computation of points  of the $C(\alpha,Q)$ curves (except circular arcs) we need:
\begin{enumerate}
\item{take an arbitrary point $\xi,\ 0\le \xi \le \pi/2$ and the tangent line $t_1$ to $Q$ in the point $T_1=(\xi, -\sin \xi)$,}
\item{compute the associated tangent line $t_2$ to $Q$ such that $\angle(t_1,t_2)=\alpha$,}
\item{the intersection point $T=(x(\xi),y(\xi))$ of tangent lines $t_1, t_2$ lies on the curve $C(\alpha,Q)$.}
\end{enumerate}

The behaviour of curve $C(\alpha,Q)$ depends on the angle $\alpha$ and it is very interesting. For $\alpha$ near (but less than) straight angle, those curves are convex. For certain obtuse angles close to $\pi/2$ the central part of curves (i.e., the part contained in the region (I)) have waveform shape. Minimal amplitudes of  ``wave oscillations'' we observed for the angle of about 1.9 radians. For obtuse angles, smaller than (approximately) 1.88 radians the above-mentioned central part is concave.
\begin{figure}[ht]
    \begin{minipage}[t]{5.5cm}
        \includegraphics[width=5cm]{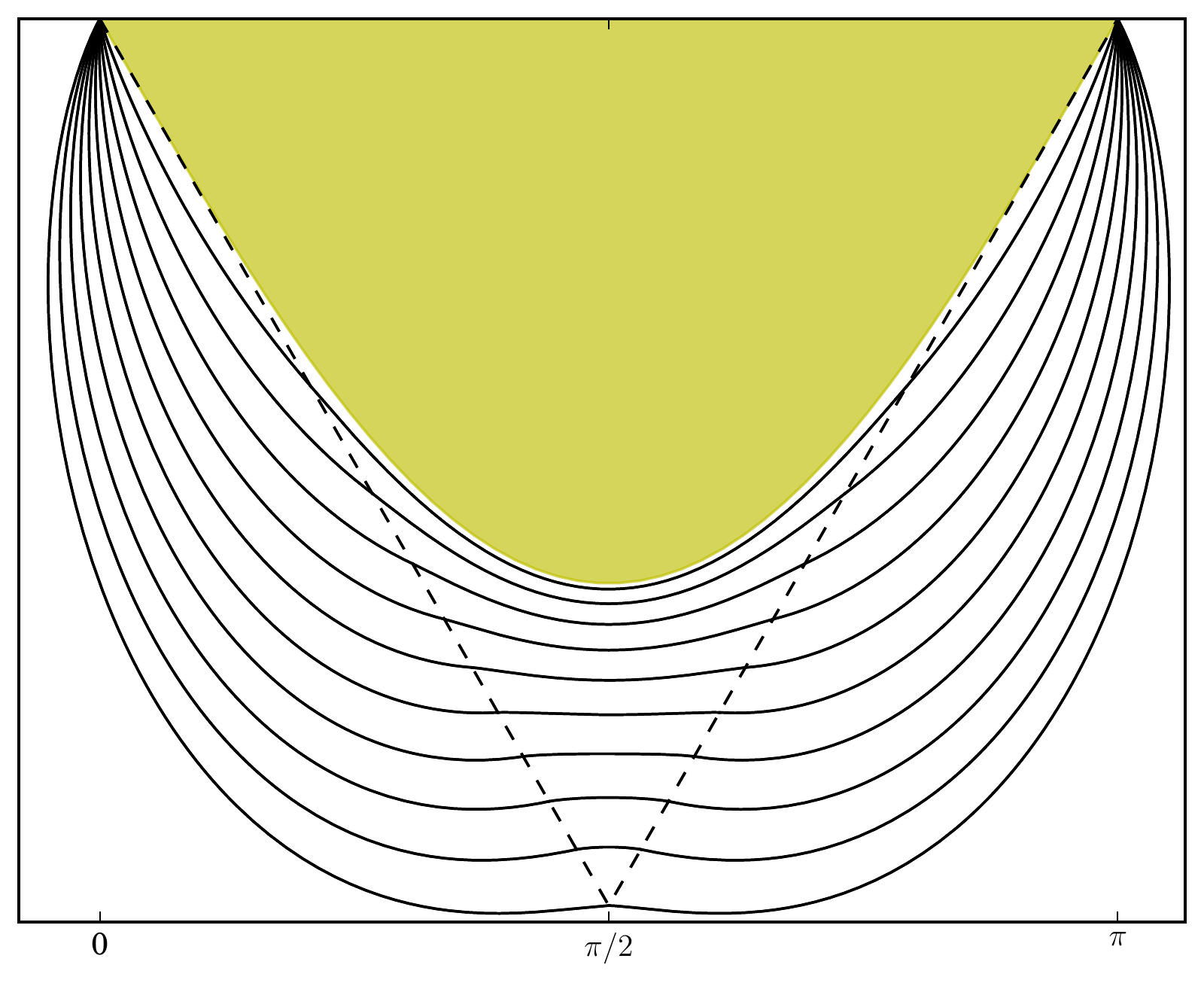}
        \caption{$C(\alpha,Q)$ for obtuse angles}
        \label{sinpicobtuse}
    \end{minipage}
     \quad
    \begin{minipage}[t]{5.5cm}
       \includegraphics[width=5cm]{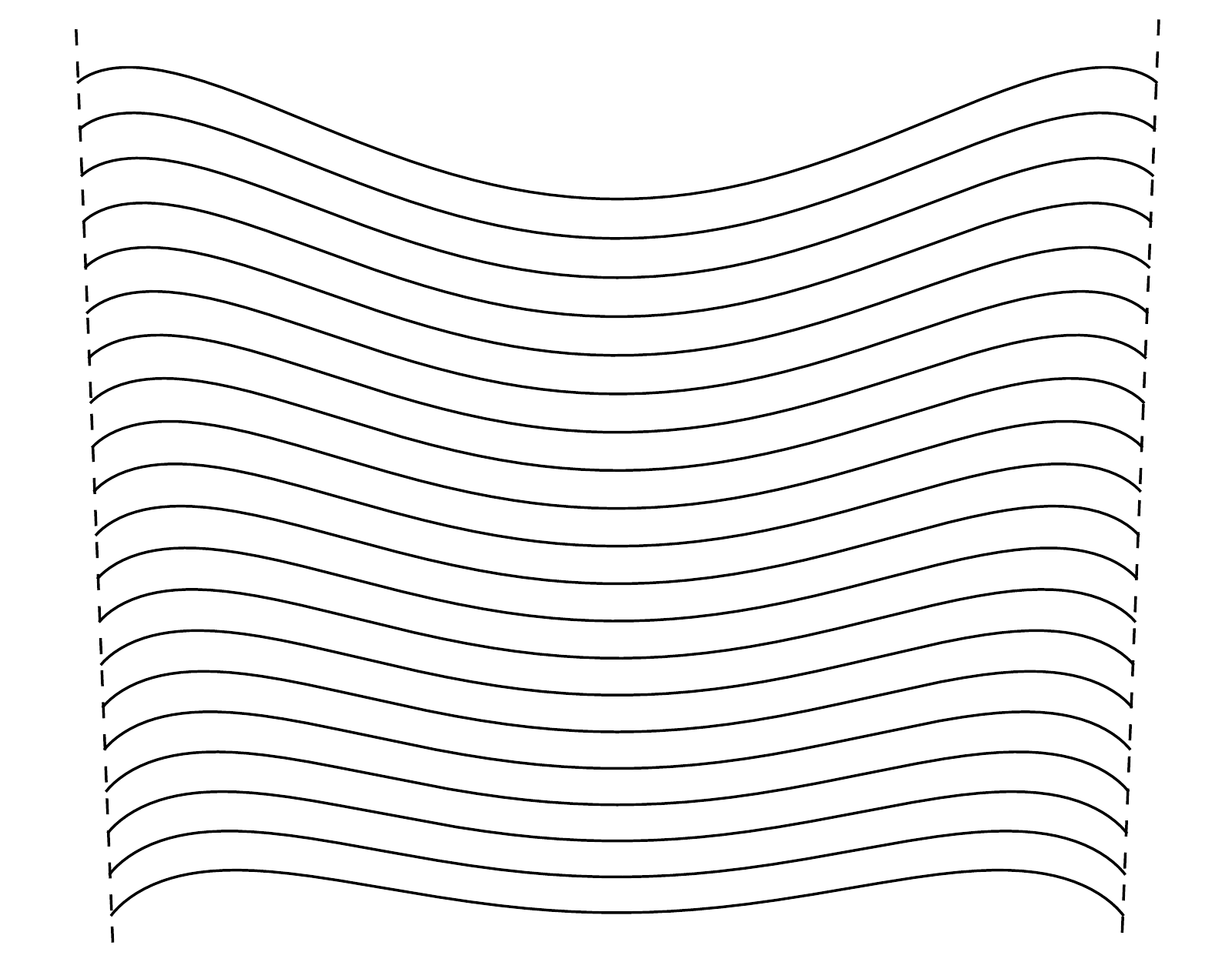}
       \caption{Wawe-shaped central parts}
       \label{sinpicwave}
       \end{minipage}
\end{figure}
On the figure \ref{sinpicobtuse} we can see the curves $C(\alpha,Q)$ for ten angles
$\alpha_k=\pi/(1+0.1\, k),\ k=1,2,\dots, 10$. Figure \ref{sinpicwave} shows central waveform parts of twenty curves  $C(\alpha,Q)$ for angles in interval from 1.92 to 1.95 radians (the picture is zoomed vertically). For sharp angles the curves $C(\alpha,Q)$ are not very interesting,
for $\alpha \rightarrow 0$ they are close to big circular arcs with central angle 
$2 \alpha$.

\begin{example}
We take the convex region $Q$, bounded by ellipse 
$$x=a \cos\varphi,\quad y=b\, \sin\varphi,\quad 0\le \varphi <2 \pi,\quad a\ge b.$$
Here, using the Sage system, we can derive the formula for $C(\alpha,Q)$ in polar coordinates. 
\begin{figure}[ht]
    \begin{minipage}[t]{5.5cm}
        \includegraphics[width=5cm]{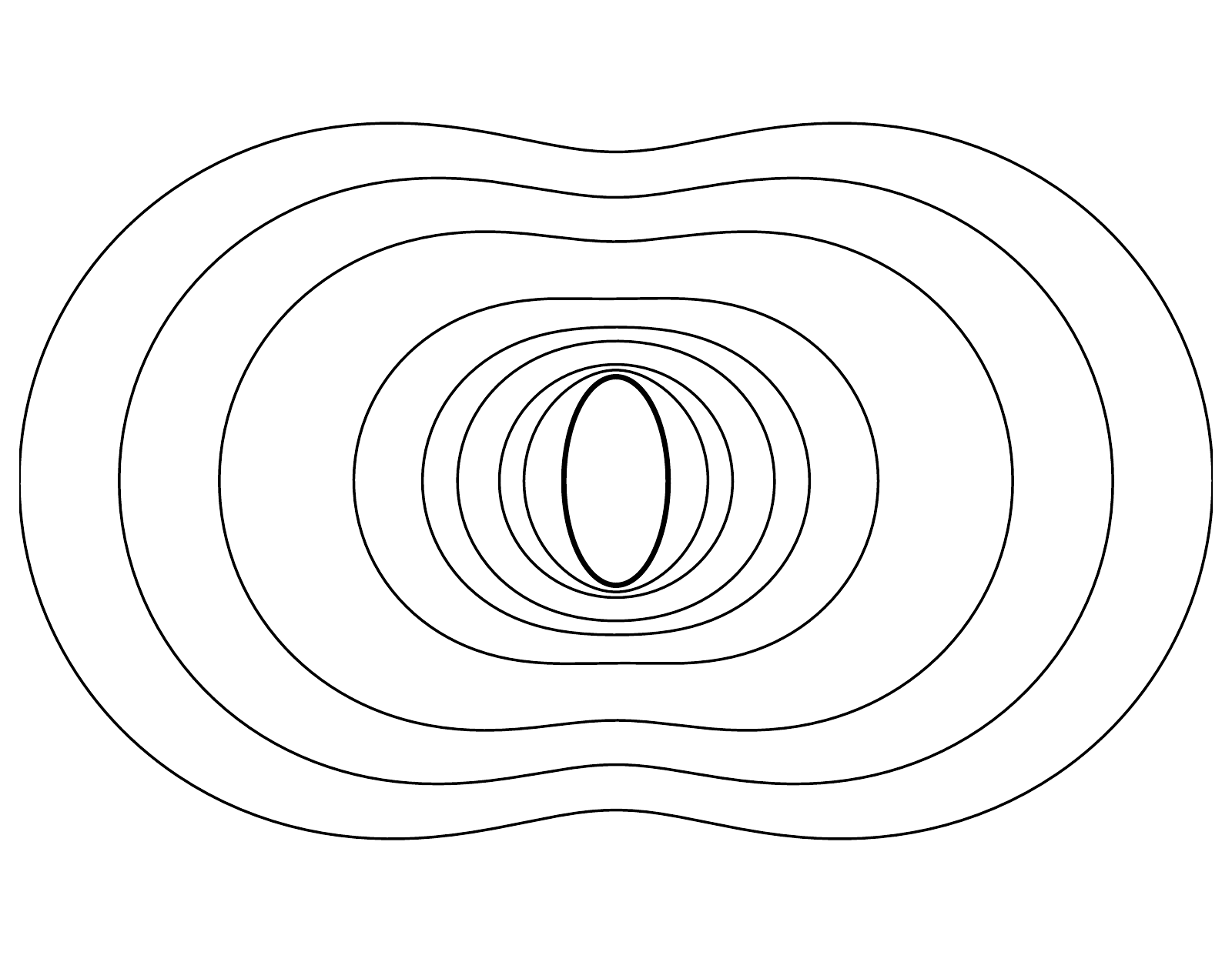}
        \caption{$C(\alpha,Q)$ for ``thin'' ellipses}
        \label{thinellipse}
    \end{minipage}
     \quad
    \begin{minipage}[t]{5.5cm}
       \includegraphics[width=5cm]{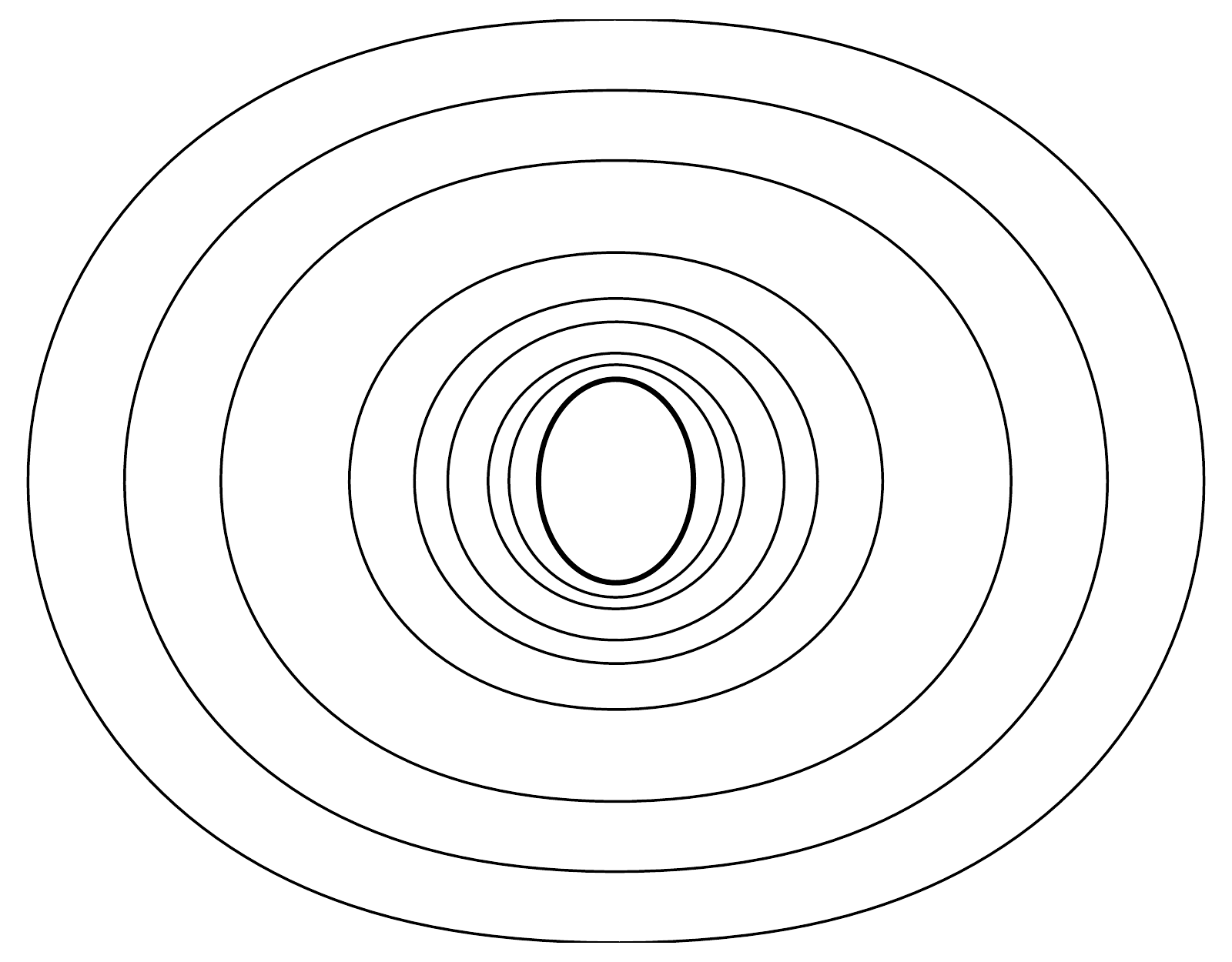}
       \caption{$C(\alpha,Q)$ for ``thick'' ellipses}
       \label{thickellipse}
       \end{minipage}
\end{figure}
Note, that the curve $C(\alpha,Q)$ for $\alpha=\pi/2$ is simply the circle with radius $r=\sqrt{a^2+b^2}$. Obviously, for angles close to straight angle, the constant view angle curves are convex. If the ratio $b/a$ is greater than approximately 0.707 (hypothesis: maybe this value is $\sqrt{2}/2$), then the regions bounded by curves $C(\alpha,Q)$ are convex for all $\alpha$. On the other side, if the ratio $b/a$ is small, the curves $C(\alpha,Q)$ for small angles will approach the two circular arcs with central angle $2 \alpha$.
\end{example}

\section*{Directions for further research}
In this paper, we brought only simple examples, but already we have seen some interesting properties of curves of equal aperture angle. The broad research area is to find efficient numerical algorithms for computations of those curves. It would be interesting to approximate $C(\alpha,Q)$ with the corresponding curve for polygonal approximation of boundary of region $Q$. Here we can explore some results from \cite{angleapprox} or \cite{teichmann} (e.g., the method of rotating wedges).


\end{document}